\nonstopmode
\pdfoutput=1
\documentclass[graybox]{svmult}

\usepackage{mathptmx}       
\usepackage{helvet}         
\usepackage{courier}        
\usepackage{type1cm}        
\usepackage{makeidx}         
\usepackage{graphicx}        
\usepackage{multicol}        
\usepackage[bottom]{footmisc}

\usepackage{amsmath,amssymb,amsfonts}
\usepackage{hyperref}
\usepackage[english]{babel}
\usepackage{todonotes}
\usepackage{enumerate}

\newcommand{\todomartin}[1]{}

\def\al{\alpha}

\def\ga{\gamma}

\def\ep{\varepsilon}

\def\th{\theta}

\def\ka{\kappa}
\def\la{\lambda}

\def\si{\sigma}

\def\ph{\varphi}

\def\ps{\psi}

\def\De{\Delta}

\def\Ph{\Phi}

\def\o{\circ}
               
\def\inv{^{-1}}
\def\x{\times}
\def\p{\partial}

\def\R{{\mathbb R}}

\def\exp{\operatorname{exp}}

\let\on=\operatorname

\let\mc=\mathcal

\newcommand{\ud}{\,\mathrm{d}}
\spnewtheorem*{OQ}{Open Question}{\bfseries}{\rmfamily}

\makeindex             

\begin{document}

\title*{Why Use Sobolev Metrics on the Space of Curves}
\author{Martin Bauer, Martins Bruveris and Peter W. Michor}
\institute{Martin Bauer,  \at University of Vienna, Faculty for Mathematics, 
Oskar-Morgenstern-Platz~1, 1090 Vienna, AT, \email{bauer.martin@univie.ac.at}
\and Martins Bruveris \at Brunel Unversity London, College of Engineering Design and Physical 
Sciences, Department of Mathematics, John Crank Building, Brunel University London, Uxbridge, UB8 
3PH, UK,  \email{martins.bruveris@brunel.ac.uk }\and  
Peter W. Michor,  \at University of Vienna, Faculty for Mathematics, Oskar-Morgenstern-Platz~1, 
1090 Vienna, AT, \email{peter.michor@univie.ac.at} 
}
\maketitle

\abstract{In this chapter we study reparametrization invariant Sobolev metrics on spaces of regular 
curves. We discuss their completeness properties and  
the resulting usability for applications in shape analysis. 
In particular, we will argue, that the development of efficient numerical methods for higher order 
Sobolev type metrics is an extremely desirable goal.}  

\section{Introduction}
Over the past decade Riemannian geometry on the infinite-dimensional spaces of parametrized and unparametrized curves has developed into an active research area. The interest has been fueled by the important role of these spaces in the areas of shape analysis and  computer vision. Contributions in these fields include applications to medical image diagnostics \cite{Gla2008},
target and activity recognition in image analysis \cite{SKKS2014}, plant leaf classification \cite{LKS2014} and protein structure analysis \cite{LSZ2011} or human motion analysis in computer graphics \cite{Esl2014}. In these research areas one is interested in studying the variability within a certain class of shapes. As a consequence an important goal is the development of statistical tools for these spaces.

Riemannian metrics provide the additional structure, that is needed to capture the nonlinearity of the space and at the same time linearize it sufficiently to enable computations. In this chapter we want to acquaint the reader with some of the metrics, that can be defined on the space of curves and discuss their properties with a view towards applications in shape analysis. We will concentrate particularly on completeness properties: do geodesics exist, when do they stop existing and how does it depend on the metric. For a more wide-ranging overview of Riemannian metrics on spaces of functions, see~\cite{Bauer2014}.

\runinhead{Parametrized Curves} In this chapter we will discuss Riemannian metrics on two different spaces: first, the space of smooth, regular, closed curves in $\R^n$
\begin{align}
\on{Imm}(S^1,\mathbb R^d)=\left\{c\in C^{\infty}(S^1,\mathbb R^d): c'(\th) \neq 0,\, \forall \th \in S^1 \right\}\;; 
\end{align}
here $\on{Imm}$ stands for \emph{immersion}. This is an open set in the Fr\'echet space $C^\infty(S^1,\R^d)$ of all smooth functions and as such it is itself a Fr\'echet manifold. As an open subset of a vector space, its tangent space at any curve is the vector space itself, 
$T\on{Imm}(S^1,\R^d) \cong \on{Imm}(S^1,\R^d) \x C^\infty(S^1,\R^d)$. A Riemannian metric on $\on{Imm}(S^1,\R^d)$ is a smooth map
\[
G : \on{Imm}(S^1,\R^d) \x C^\infty(S^1,\R^d) \x C^\infty(S^1,\R^d) \to \R\;,
\]
such that $G_c(\cdot, \cdot)$ is a symmetric, positive definite bilinear form for all curves $c$. An example of a Riemannian metric is the $L^2$-metric $G_c(h,k) = \int_{S^1} \langle h, k \rangle |c'| \ud \th$, which we will look at more closely in Sect.~\ref{sec:L2metric}. When studying particular Riemannian metrics, it will be useful to consider larger spaces of less regular curves, but $\on{Imm}(S^1,\R^d)$ will always be the common core.

\runinhead{Unparametrized Curves} The other space, that we will consider, is the space of unparametrized curves, sometimes also denoted \emph{shape space}. There are several closely related, but slightly differing ways to define this space mathematically. We will consider an \emph{un\-pa\-ra\-met\-rized curve} or \emph{shape} to be an equivalence class of parametrized curves, that differ only by a reparametrization. In other words, $c_1$ and $c_2$ represent the same shape, if $c_1 = c_2 \o \ph$ for some reparametrization $\ph \in \on{Diff}(S^1)$; mathematically $\on{Diff}(S^1)$ is the diffeomorphism group of the circle, that is, the set of all smooth invertible maps $\ph : S^1\to S^1$. With this definition the space of unparametrized curves is the quotient
\[
B(S^1,\R^d) = \on{Imm}(S^1,\R^d) / \on{Diff}(S^1)\;.
\]
Apart from isolated singular points, the space $B(S^1,\R^d)$ is also an infinite-di\-men\-sional manifold and the projection $p: c \to [c]$ assigning each curve its equivalence class is a submersion.\footnote{In applications one  often wants to consider curves and shapes modulo Euclidean motions, leading to the spaces $\on{Imm}(S^1,\R^d)/\on{Mot}$ and $B(S^1,\R^d)/\on{Mot}$, where $\on{Mot} = SO(d) \ltimes \R^d$ denotes the Euclidean motion group. All metrics discussed in this chapter are invariant under the motion group and therefore induce a Riemannian metric on the quotients $\on{Imm}(S^1,\R^d)/\on{Mot}$ and $B(S^1,\R^d)/\on{Mot}$. In Sect.~\ref{sec:SRVT} we will encounter metrics, that live naturally on the space $\on{Imm}(S^1,\R^d)/\on{Tra}$ of curves modulo translations.}

\runinhead{Reparametrization Invariant Metrics}
To define a Riemannian metric on the space $B(S^1,\R^d)$ we will start with a Riemannian metric $G$ on $\on{Imm}(S^1,\R^d)$, that is invariant under the action of $\on{Diff}(S^1)$; such metrics are called \emph{reparametrization invariant}. This means, $G$ has to satisfy
\[
G_{c\o \ph}(h \o \ph, k \o \ph) = G_c(h,k)\;,
\]
for all curves $c$, tangent vectors $h,k$ and reparametrizations $\ph$. Then we can use the formula
\[
G_{[c]}(X,X) = \inf \left\{ G_c(h,h) \,:\, T_c p.h = X \right\}\;,
\]
to define a Riemannian metric on shape space $B(S^1,\R^d)$, such that the projection $p$ is a Riemannian submersion.

\runinhead{Geodesic Distance}
An important concept in shape analysis is the notion of distance between two curves or shapes. A Riemannian metric leads to a natural distance function, the \emph{induced geodesic distance}. The distance measures the length of the shortest path between two curves. If $c_0, c_1 \in \on{Imm}(S^1,\R^d)$ are two parametrized curves, the distance between them is defined as
\[
\on{dist}_I(c_0, c_1) = \inf_{\substack{\ga(0) = c_0 \\ \ga(1) = c_1}} \int_0^1 \sqrt{G_{\ga(t)}(\ga_t(t), \ga_t(t))} \ud t\;,
\]
where the infimum is taken over all smooth paths $\gamma$, that connect the curves $c_0$ and $c_1$. Whether there exists a path, realizing this infimum, is an interesting and non-trivial question in Riemannian geometry.

If we start with a reparametrization invariant metric $G$ and it induces a Riemannian metric on shape space $B(S^1,\R^d)$, then we will be interested in computing the geodesic distance on $B(S^1,\R^d)$. The geodesic distances for parametrized and unparametrized curves are related by
\begin{equation}
\on{dist}_B([c_0],[c_1]) = \underset{\varphi \in \on{Diff}(S^1)}{\on{inf}} \on{dist}_I(c_0,c_1\circ\varphi)\;.
\end{equation}
For applications in shape analysis it is important to find a stable and fast method 
to numerically compute this quantity for arbitrary shapes $[c_0]$ and $[c_1]$. We will comment in the later sections, for which metrics a reparametrization $\ph$, realizing the above infimum, exists, and what the obstructions to its existence are in other cases.

\runinhead{Organization} We will look at three families of metrics: first, the $L^2$-metric in Sect.~\ref{sec:L2metric}, which is among the simplest reparametrization invariant metrics, but unfortunately unsuitable for shape analysis; then, first order Sobolev metrics in Sect.~\ref{sec:SRVT}, which are very well suited for numerical computations and therefore among the most widely used Riemannian metrics in applications; finally, we will look at higher order Sobolev metrics in Sect.~\ref{sec:higherorder} and argue, why their theoretical properties make them good candidates for use in shape analysis. At the end we will explain, how these metrics can be generalized to spaces of parametrized and unparametrized surfaces.

\section{The $L^2$-Metric}\label{sec:L2metric}

The arguably simplest Riemannian metric on the space of smooth, regular curves, that is invariant under reparametrizations, is the $L^2$-metric
\[
G_c(h,k)=\int_{S^1}\langle h,k\rangle \ud s\;,
\]
where we use $\ud s = |c'| \ud \th$ to denote arc length integration. It is integration with respect to $\ud s$ rather than $\ud \th$, that makes this metric reparametrization invariant, as can be seen from the following calculation,
\begin{equation*}
 G_{c\circ \ph}(h\circ\ph,k\circ\ph) 
= \int_{S^1} \langle h\circ\ph, k\circ\ph \rangle \, (|c'|\circ\ph) \, \ph' \ud \th
= G_c(h,k)\;.
\end{equation*}
Similarly, if we wanted to include derivatives of $h, k$ in the metric and keep the metric reparametrization invariant, we would need to use the arc length differentiation $D_s h = \frac1{|c'|} h'$ rather that $h' = \p_\th h$.

\runinhead{Geodesic Equation}
The geodesic equation of the $L^2$-metric is a nonlinear, second order PDE for the path $c(t,\th)$. It has the form
\begin{equation}\label{eq:geode_L2}
\left(|c_{\theta}|c_t\right)_t = -\frac12 \left(\frac{|c_t|^2}{|c_{\theta}|}\, c_{\theta}\right)_{\theta}\;.
\end{equation}
where $c_\th = \p_\th c = c'$ and $c_t = \p_t c$ denote the partial derivatives. While the equation is as simple as one can hope for---the geodesic equations for higher order metrics have many more terms---there are currently no existence results available for it.
\begin{OQ}
Given a pair of an initial curve and an initial velocity $(c_0, u_0) \in T\on{Imm}(S^1,\R^d)$, does the geodesic equation admit short time solutions with the given initial conditions?
\end{OQ}
We know, that we cannot hope for long time existence, since it is possible to shrink a circle along a geodesic path down to a point in finite time. Numerical evidence in \cite[Sect.~5.3] {Michor2006c} suggest that geodesics should exist as long as the curvature of the curve remains bounded.

\runinhead{Geodesic Distance}
The lack of existence results for the geodesic equation is not the biggest problem of the $L^2$-metric. The crucial property, that makes it unsuitable for applications in shape analysis, is that the induced geodesic distance vanishes. 

The \emph{geodesic distance} between two curves $c_0, c_1 \in \on{Imm}(S^1,\R^d)$ is defined as the infimum over the lengths of all paths $\ga$, connecting the two curves, i.e.
\[
\on{dist}_{I}(c_0, c_1) = \inf_{\substack{\ga(0) = c_0 \\ \ga(1) = c_1}} \int_0^1 \sqrt{G_{\ga(t)}(\ga_t(t), \ga_t(t))} \ud t\;.
\]
It was found \cite{Michor2006c,Michor2005,Bauer2012c} that for the $L^2$-metric the geodesic 
distance between any two curves is zero.\footnote{We encounter the vanishing of the geodesic 
distance for $L^2$-metrics on several spaces: on the space $\on{Imm}(M,N)$ of immersions between 
two manifolds $M$, $N$ of arbitrary dimension, $M$ compact; on the Virasoro--Bott group 
\cite{Bauer2012c}; and even for Sobolev metrics on the diffeomorphism group of a compact manifold, 
provided the order of the metric is $< \frac 12$ \cite{Bauer2013b, Bauer2013c}.} What does this 
mean? If $\ga$ is a smooth, non-constant path, then $\p_t\ga(t)$ cannot be identically zero and so 
the length $\int_0^1 \sqrt{G_\ga(\ga_t, \ga_t)} \ud t$ will be strictly positive. The meaning of 
$\on{dist}_I(c_0, c_1) =0$ is that we can find arbitrary short paths connecting $c_0$ and $c_1$. No 
path will have zero length, but given any $\ep > 0$, we can find a path with length $< \ep$. How do 
these paths look like? They are easier to visualize for the geodesic distance on the space of 
unparametrized curves, which we will describe next. 

For the $L^2$-metric the geodesic distance between the unparametrized curves $[c_0], [c_1] \in B(S^1,\R^d)$, represented by $c_0, c_1$ can be computed as the following infimum,
\[
\on{dist}_B([c_0], [c_1]) = \inf_\ga \int_0^1 \sqrt{ G_\ga(\ga_t^\perp, \ga_t^\perp)} \ud t\;;
\]
here $\ga(t)$ is a path starting at $c_0$ and ending at any curve in the equivalence class $[c_1]$, that is $\ga(1) = c_1 \o \ph$ for some $\ph \in \on{Diff}(S^1)$. We denote by $\ga_t^\perp = \ga_t - \langle \ga_t, v \rangle v$, with $v = D_s \ga$, the projection of the vector $\ga_t(t,\th) \in \R^d$ to the subspace orthogonal to the curve $\ga(t)$ at $\th$.

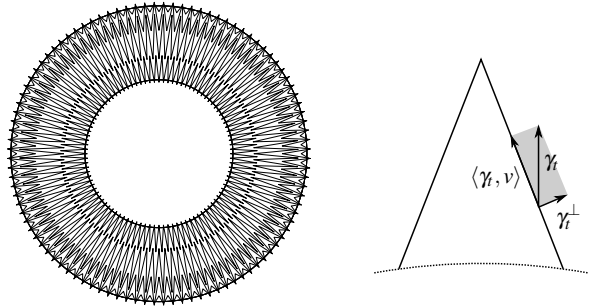
\begin{figure}
\begin{center}
\includegraphics[width=.35\textwidth]{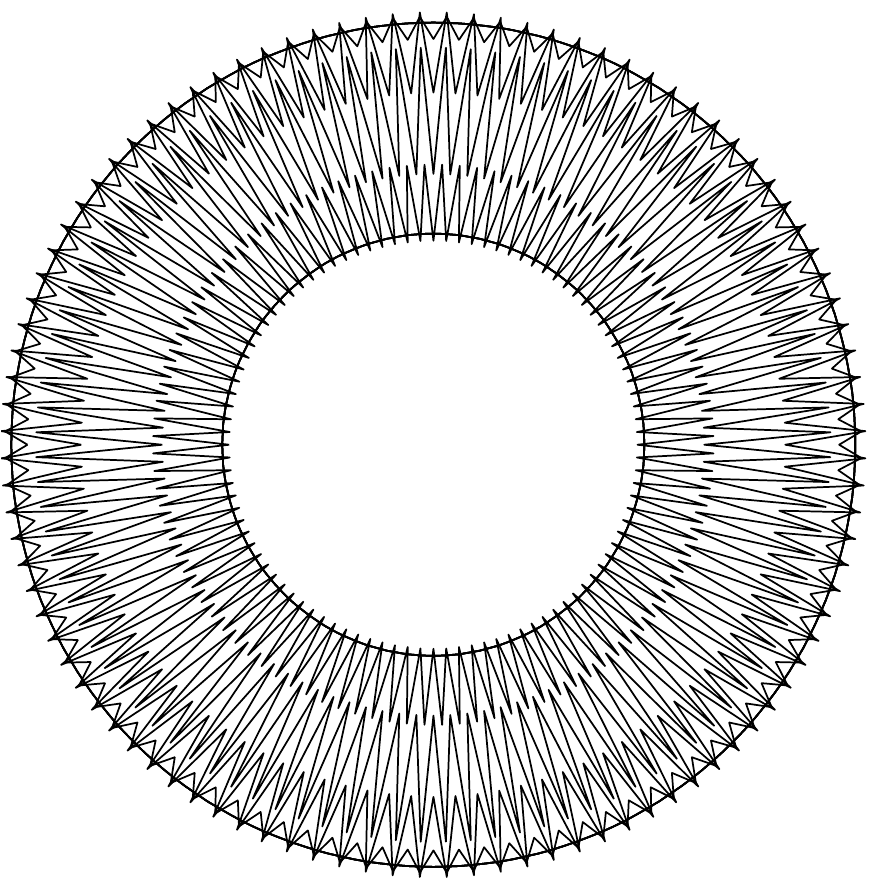}
\def\svgwidth{.35\columnwidth}
\begingroup%
  \makeatletter%
  \providecommand\color[2][]{%
    \errmessage{(Inkscape) Color is used for the text in Inkscape, but the package 'color.sty' is not loaded}%
    \renewcommand\color[2][]{}%
  }%
  \providecommand\transparent[1]{%
    \errmessage{(Inkscape) Transparency is used (non-zero) for the text in Inkscape, but the package 'transparent.sty' is not loaded}%
    \renewcommand\transparent[1]{}%
  }%
  \providecommand\rotatebox[2]{#2}%
  \ifx\svgwidth\undefined%
    \setlength{\unitlength}{300bp}%
    \ifx\svgscale\undefined%
      \relax%
    \else%
      \setlength{\unitlength}{\unitlength * \real{\svgscale}}%
    \fi%
  \else%
    \setlength{\unitlength}{\svgwidth}%
  \fi%
  \global\let\svgwidth\undefined%
  \global\let\svgscale\undefined%
  \makeatother%
  \begin{picture}(1,1)%
    \put(0,0){\includegraphics[width=\unitlength]{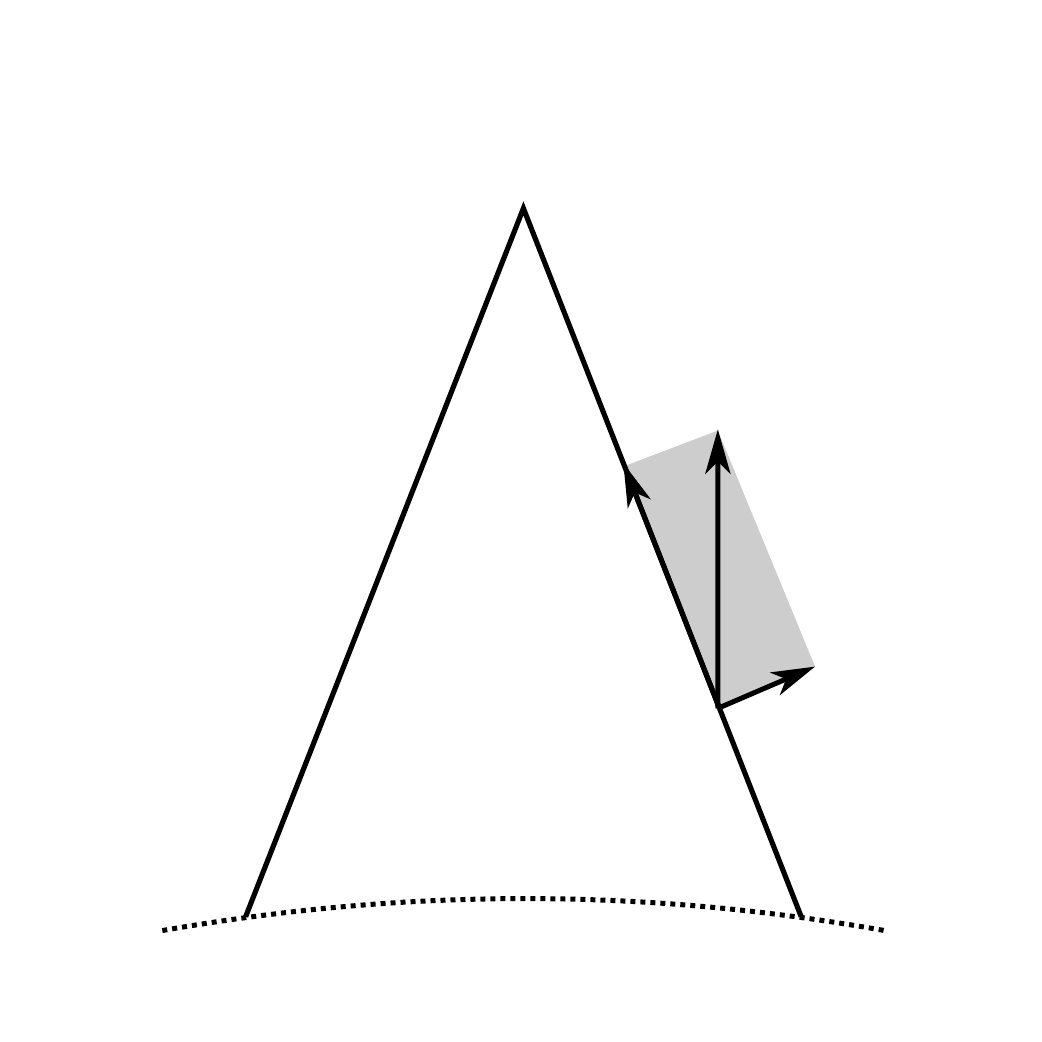}}%
    \put(0.74666667,0.26666667){\color[rgb]{0,0,0}\makebox(0,0)[lb]{\smash{$\gamma_t^\bot$}}}%
    \put(0.46666667,0.4){\color[rgb]{0,0,0}\makebox(0,0)[lb]{\smash{$\langle \gamma_t, v \rangle$}}}%
    \put(0.70666667,0.45333333){\color[rgb]{0,0,0}\makebox(0,0)[lb]{\smash{$\gamma_t$}}}%
  \end{picture}%
\endgroup%

\end{center}
\caption{Left side: a short curve with respect to the $L^2$-metric in the space $B(S^1,\R^2)$ of unparametrized curves connecting two concentric circles. We see that the intermediate curves are sawtooth-shaped. Right side: Along a sawtooth the tangential component $\langle \ga_t, v \rangle v$ is large, while the normal component $\ga^\perp$ becomes small, the steeper the slope of the sawtooth.
}
\label{fig:L2metric}
\end{figure}

A short path connecting two concentric circles can be seen in Fig.~\ref{fig:L2metric}. The key observation is that the sawtooth-shaped curves have a large tangential velocity, but only a small normal velocity. Since for the geodesic distance on $B(S^1,\mathbb R^d)$ we only measure the normal part of the velocity vector, these paths have a short length. The more teeth we use, the smaller the normal component of the velocity and the smaller the length of these paths.

The vanishing of the geodesic distance started the search for stronger metrics, that would be more useful to shape analysis.

\runinhead{Almost Local Metrics}
One class of metrics, designed to have non-vanishing distance, while being as simple as possible, is the class of almost local metrics.
The motivating idea behind almost local metrics was the observation that for paths with short length in the $L^2$-metric, the intermediate curves are long and have large curvature.
Thus one hopes that by adding weights, that depend on length and curvature, to the metric, these paths will be sufficiently penalized and the geodesic distance will become non-zero. \emph{Almost local metrics}\footnote{These metrics are not local, because the length $\ell_c$ is not a local quantity, however it is only a mild non-locality; hence the name ``almost local'' metrics.} are metrics of the form
\begin{equation}
G_c(h,k)=\int_{S^1}\Ph(\ell_c,\ka)\langle h,k\rangle \ud s\;,
\end{equation}
with $\Ph$ some function of the two variables $\ell_c=\int_{S^1}\ud s$ (length) and 
$\ka$ (curvature).
If $\Ph$ depends only on $\ell_c$, the resulting metric
$G_c(h,k) = \Ph(\ell_c) \int_{S^1} \langle h, k \rangle \ud s$ is a conformal rescaling of the $L^2$-metric \cite{YezziMennucci2005,Shah2008}. Other choices for $\Ph$ include $\Ph(\ka) = 1 + A\ka^2$ with $A$ a positive constant \cite{Michor2006c} or 
the scale invariant metric $\Ph(\ell_c,\ka) = \frac{1}{\ell^3_c}+{\ka^2}{\ell_c}$ \cite{Michor2007}. 

For all these metrics it has been shown that they induce a point-separating distance
function\footnote{A distance function $d(\cdot,\cdot)$ is {\itshape point-separating}, if $d(x,y) > 0$ whenever $x \neq y$. This is stronger than non-vanishing, which would only require two points $x,y$ with $d(x,y) \neq 0$.} 
on the space $B(S^1,\mathbb R^d)$ of unparametrized curves. However, similarly to the $L^2$-metric, little is known about solutions of the geodesic equation and while the geodesic distance is point-separating on the space $B(S^1,\R^d)$, it is not point-separating on the space $\on{Imm}(S^1,\R^d)$ of parametrized curves. In the next two sections we will discuss a different strategy to strengthen the $L^2$-metric, leading to the class of Sobolev metrics.

\section{First Order Metrics and the Square Root Velocity Transform}\label{sec:SRVT}

One way to deal with the degeneracy of the $L^2$-metric is by adding terms, that involve first derivatives of the tangent vectors. Such metrics are called \emph{first order Sobolev metrics} or, short, $H^1$-metrics. An example is the metric
\[
G_c(h, k) = \int_{S^1} \langle h, k \rangle + \langle D_s h, D_s k \rangle \ud s\;,
\]
with $D_s h = \frac 1{|c'|} h'$ denoting the arc length derivative and $\ud s = |c'| \ud \th$. If we omit the $L^2$-term, we arrive at 
$G_c(h, k) = \int_{S^1} \langle D_s h, D_s k \rangle \ud s$, 
which is a metric on the space $\on{Imm}(S^1,\R^d)/\on{Tra}$ of regular curves modulo translations. The scale-invariant version of this metric has been studied in \cite{Younes1998, Michor2008a} and it has the remarkable property that one can find explicit formulas for minimizing geodesics between any two curves.

We will concentrate in this section on a related metric, obtained by using different weights for the tangential and normal components of $D_s h$,
\begin{equation}
\label{eq:srvt_metric}
G_c(h, k) = \int_{S^1} \langle D_s h^\perp, D_s k^\perp \rangle 
+ \frac 14 \langle D_s h,v\rangle \langle D_s k,v\rangle \ud s\;;
\end{equation}
here $v = D_s c = \frac 1{|c'|}c'$ is the unit length tangent vector along $c$ and 
$D_s h^\perp = D_s h - \langle D_s h, v \rangle v$ is the projection of $D_s h$ to the subspace $\{ v \}^\perp$ orthogonal to the curve. This is a Riemannian metric on $\on{Imm}(S^1,\R^d)/\on{Tra}$ and it is the metric used in the \emph{square root velocity (SRV)} framework \cite{Jermyn2011}. The reason for singling out this metric is that the SRV framework has been used successfully in applications \cite{SKKS2014,XKS2014,LKS2014} and the SRV transform has a simple and accessible form. We will comment on other $H^1$-metrics at the end of the section.

\runinhead{The Square Root Velocity Transform}
The \emph{square root velocity transform (SRVT)}  is the map
\begin{align*}
R: \operatorname{Imm}(S^1,\mathbb R^d) \rightarrow C^{\infty}(S^1,\mathbb R^d)\,,\quad
c \mapsto \frac{1}{\sqrt{|c'|}} c'\;,
\end{align*}
assigning each curve $c$ a function $q = R(c)$. 

Every vector space $(V, \langle \cdot, \cdot \rangle)$ with an inner product can be regarded as a Riemannian manifold: the Riemannian metric $g$ at each point $x \in V$ is simply the inner product, $g_x(\cdot, \cdot) = \langle \cdot,\cdot \rangle$. The SRVT is an isometry between the Riemannian manifold $\left(\on{Imm}(S^1,\R^d)/\on{Tra}, G\right)$, where $G$ is the Riemannian metric~\eqref{eq:srvt_metric} and the space $C^\infty(S^1,\R^d)$ with the $L^2$-inner product 
$\langle u, v \rangle_{L^2} = \int_{S^1} \langle u, v \rangle \ud \th$.

\runinhead{The SRVT for Open Curves}
Things are simple on the space of open curves $\on{Imm}([0,2\pi],\R^d)/\on{Tra}$. The SRVT is a one-to-one mapping between the space $\on{Imm}([0,2\pi],\R^d)/\on{Tra}$ and the set 
$C^\infty(S^1,\R^d \setminus \{0\})$ of functions that don't pass through the origin in $\R^d$. This is an open subset of all functions and thus geodesics with respect to the metric~\eqref{eq:srvt_metric} correspond to straight lines under the SRVT: the path $c(t) = R\inv(q_0 + th)$ is a geodesic in $\on{Imm}([0,2\pi], \R^2)/\on{Tra}$ and given two curves $c_0, c_1$, the geodesic connecting them is
\[
c(t) = R\inv((1-t)q_0 + tq_1)\;,
\]
with $q_i = R(c_i)$.

\runinhead{The SRVT for Closed Curves}
Things are slightly more complicated for closed curves. The inverse of the SRVT is given by the formula
\[
R\inv(q)(\th) = \int_0^\th q|q| \ud \si\;,
\]
and we see that if we want the curve $c = R\inv(q)$ to be closed, i.e., $c(0) = c(2\pi)$, then we need $\int_{S^1} q |q| \ud \si = 0$. Indeed the image of $\on{Imm}([0,2\pi],\R^d)/\on{Tra}$ under the SRVT is the set
\begin{equation*}
\on{Im}(R)= \left\{q\in C^{\infty}(S^1,\mathbb R^2):\;q(\theta)\neq 0\;\text{and}\; \int_{S^1}|q|q\;d\theta=0     \right\}\;.
\end{equation*}
We have the condition $q'(\th) \neq 0$ as before to ensure that $c'(\th) \neq 0$ and an additional constraint so that the curves $c=R\inv(q)$ are closed. Even though we don't have a closed expression for the geodesics, it is still possible to compute the geodesics numerically without much difficulty.

\runinhead{Minimizing Geodesics for Parametrized Curves}
Let us first look at open curves. In the SRV representation, $q_i = R(c_i)$, the minimizing path between $q_0$ and $q_1$ is the straight line $q(t) = (1-t)q_0 + t q_1$. In particular the minimizing path always exists. It can happen, however, that the straight line between $q_0(\th)$ and $q_1(\th)$ passes through the origin; at points $(t,\th)$, where this happens the derivative of the curve $c(t,\th) = R\inv(q(t))(\th)$ vanishes, i.e., $c'(t,\th) = 0$ and thus the curve is not regular at those points. Apart from that, any two curves can be joined by a unique minimizing geodesic, which can be computed via an explicit formula, and we know when the intermediate curves will fail to be regular.

\begin{figure}
\begin{center}
\includegraphics[width=.96\textwidth]{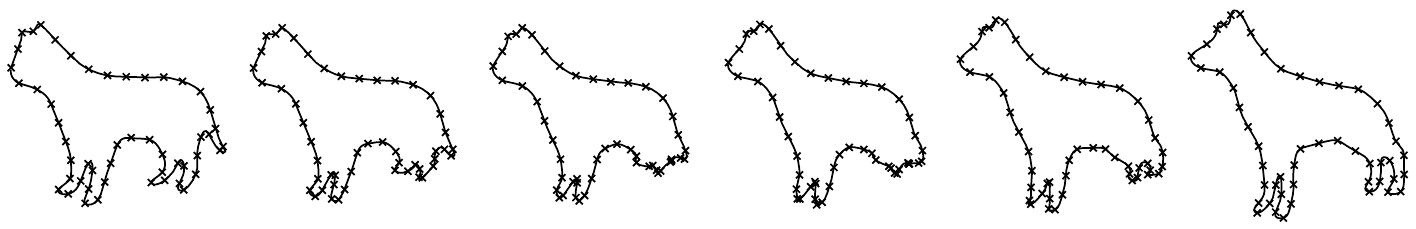}
\includegraphics[width=.96\textwidth]{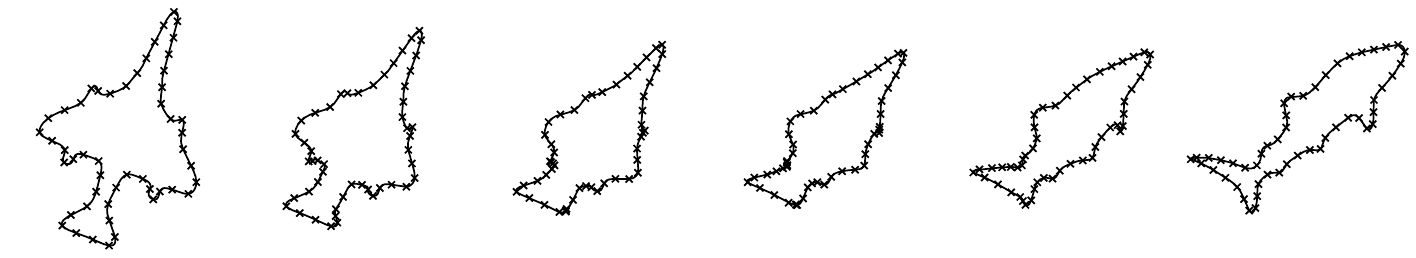}
\end{center}
\caption{Minimal geodesics between two pairs of parametrized curves. Images taken from \cite{Bauer2014b}.}
\label{fig:geod_param_curves}
\end{figure}

For closed curves the situation is less explicit, because now we also have to satisfy the nonlinear constraint $\int_{S^1} q |q| \ud \th = 0$. This is a $d$-dimensional constraint on an otherwise infinite-dimensional space and furthermore the function $q \mapsto \int_{S^1} q |q| \ud \th$ is continuous with respect to the $L^2$-topology. Numerical evidence suggest, that minimizing geodesics continue to exist between any two curves.
In particular, computing minimizing geodesics between parametrized curves is a fast and stable operation; an example of a geodesic can be seen in Fig.~\ref{fig:geod_param_curves}.

Smoothness of the minimizing geodesics is another issue. The natural target space for the SRVT is 
the space $L^2(S^1,\R^d)$ of square-integrable functions. If the SRVT of a curve lies in 
$L^2(S^1,\R^d)$, the curve itself is only absolutely continuous. Unfortunately the Riemannian 
metric $G_c(h,k)$, given by~\eqref{eq:srvt_metric}, does not have to be finite for absolutely 
continuous curves $c$ and tangent vectors $h,k$; the term $D_s h = \frac{1}{|c'|} h'$ may well 
become infinite. We are approaching the frontier of the Riemannian framework now: any two (open) curves 
can be joined by a minimizing path, however the space, where the path lives---the completion of the 
space of smooth curves, if one wants to use the term---is not a Riemannian manifold any more.  

\runinhead{Minimizing Geodesics for Unparametrized Curves}
If we want to find minimizing geodesics between two unparametrized curves $C_0, C_1 \in B(S^1,\R^d)/\on{Tra}$, represented by the curves $c_0, c_1 \in \on{Imm}(S^1,\R^d)/\on{Tra}$, one way to do this is to minimize $\on{dist}_{I}(c_0, c_1 \o \ph)$ over $\ph \in \on{Diff}(S^1)$ or equivalently over all parametrized curves $c_1 \o \ph$ representing the shape $C_1$; indeed, the geodesic distance on $B(S^1,\R^d)/\on{Tra}$ is given by
\begin{equation}
\label{eq:shape_distance}
\on{dist}_B(C_0, C_1) = \inf_{\ph \in \on{Diff}(S^1)} \on{dist}_I(c_0, c_1 \o \ph)\;.
\end{equation}
If the infimum is attained for $\ps \in \on{Diff}(S^1)$ and if we denote by $c(t)$ the minimizing geodesic between $c_0$ and $c_1 \o \ps$, then the curve $[c(t)]$ in $B(S^1,\R^d)/\on{Tra}$ is the minimizing geodesic between $C_0$ and $C_1$. Thus we are interested, whether the infimum~\eqref{eq:shape_distance} is attained and if it is, in what space.

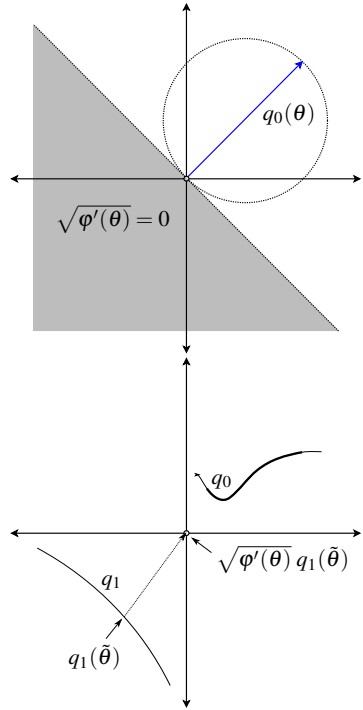
\begin{figure}
\centering
\def\svgwidth{.40\columnwidth}
\begingroup%
  \makeatletter%
  \providecommand\color[2][]{%
    \errmessage{(Inkscape) Color is used for the text in Inkscape, but the package 'color.sty' is not loaded}%
    \renewcommand\color[2][]{}%
  }%
  \providecommand\transparent[1]{%
    \errmessage{(Inkscape) Transparency is used (non-zero) for the text in Inkscape, but the package 'transparent.sty' is not loaded}%
    \renewcommand\transparent[1]{}%
  }%
  \providecommand\rotatebox[2]{#2}%
  \ifx\svgwidth\undefined%
    \setlength{\unitlength}{307bp}%
    \ifx\svgscale\undefined%
      \relax%
    \else%
      \setlength{\unitlength}{\unitlength * \real{\svgscale}}%
    \fi%
  \else%
    \setlength{\unitlength}{\svgwidth}%
  \fi%
  \global\let\svgwidth\undefined%
  \global\let\svgscale\undefined%
  \makeatother%
  \begin{picture}(1,1)%
    \put(0,0){\includegraphics[width=\unitlength]{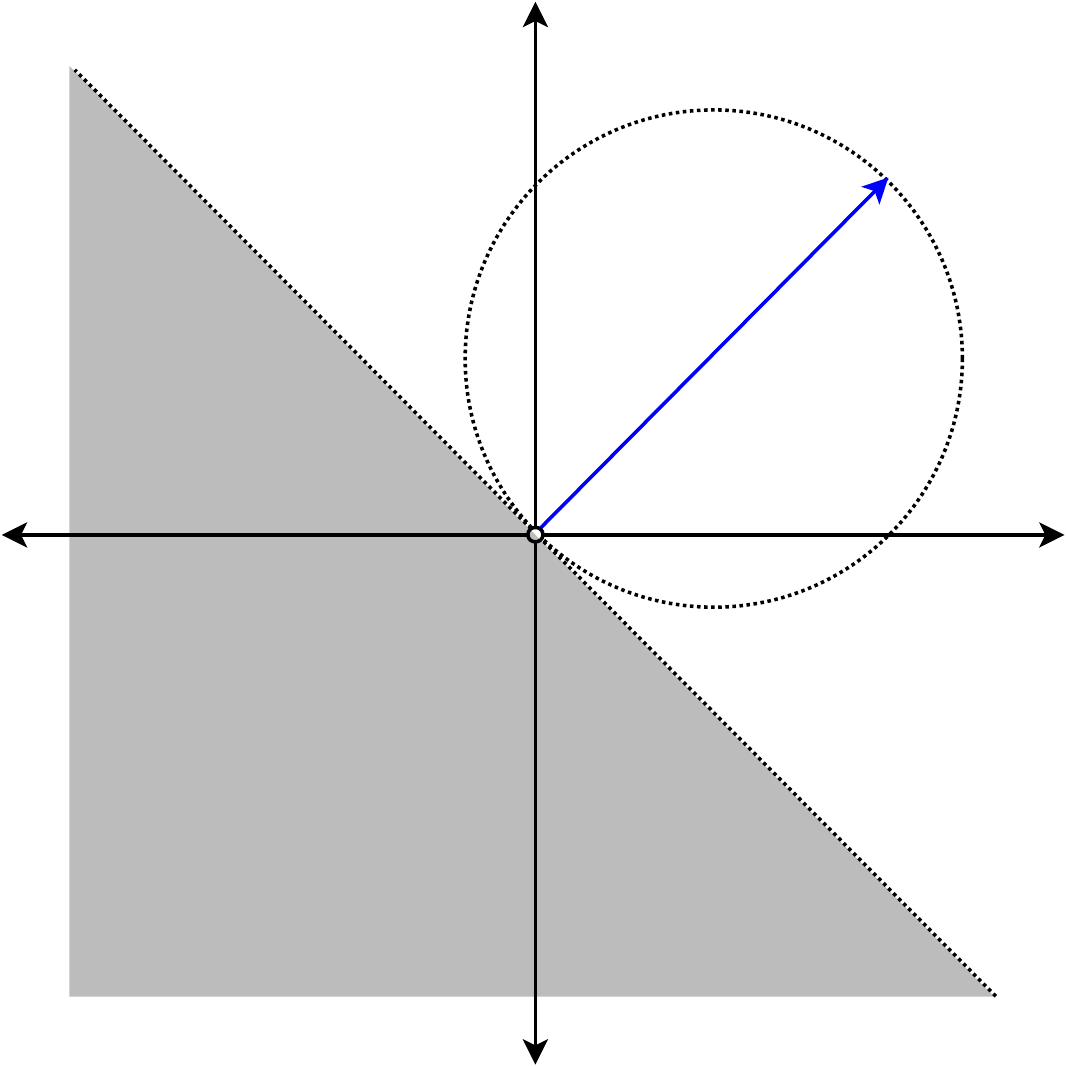}}%
    \put(0.13029316,0.36482085){\color[rgb]{0,0,0}\makebox(0,0)[lb]{\smash{$\sqrt{\ph'(\theta)}=0$}}}%
    \put(0.71921824,0.66579805){\color[rgb]{0,0,0}\makebox(0,0)[lb]{\smash{$q_0(\theta)$}}}%
  \end{picture}%
\endgroup%

\def\svgwidth{.40\columnwidth}
\begingroup%
  \makeatletter%
  \providecommand\color[2][]{%
    \errmessage{(Inkscape) Color is used for the text in Inkscape, but the package 'color.sty' is not loaded}%
    \renewcommand\color[2][]{}%
  }%
  \providecommand\transparent[1]{%
    \errmessage{(Inkscape) Transparency is used (non-zero) for the text in Inkscape, but the package 'transparent.sty' is not loaded}%
    \renewcommand\transparent[1]{}%
  }%
  \providecommand\rotatebox[2]{#2}%
  \ifx\svgwidth\undefined%
    \setlength{\unitlength}{307bp}%
    \ifx\svgscale\undefined%
      \relax%
    \else%
      \setlength{\unitlength}{\unitlength * \real{\svgscale}}%
    \fi%
  \else%
    \setlength{\unitlength}{\svgwidth}%
  \fi%
  \global\let\svgwidth\undefined%
  \global\let\svgscale\undefined%
  \makeatother%
  \begin{picture}(1,1)%
    \put(0,0){\includegraphics[width=\unitlength]{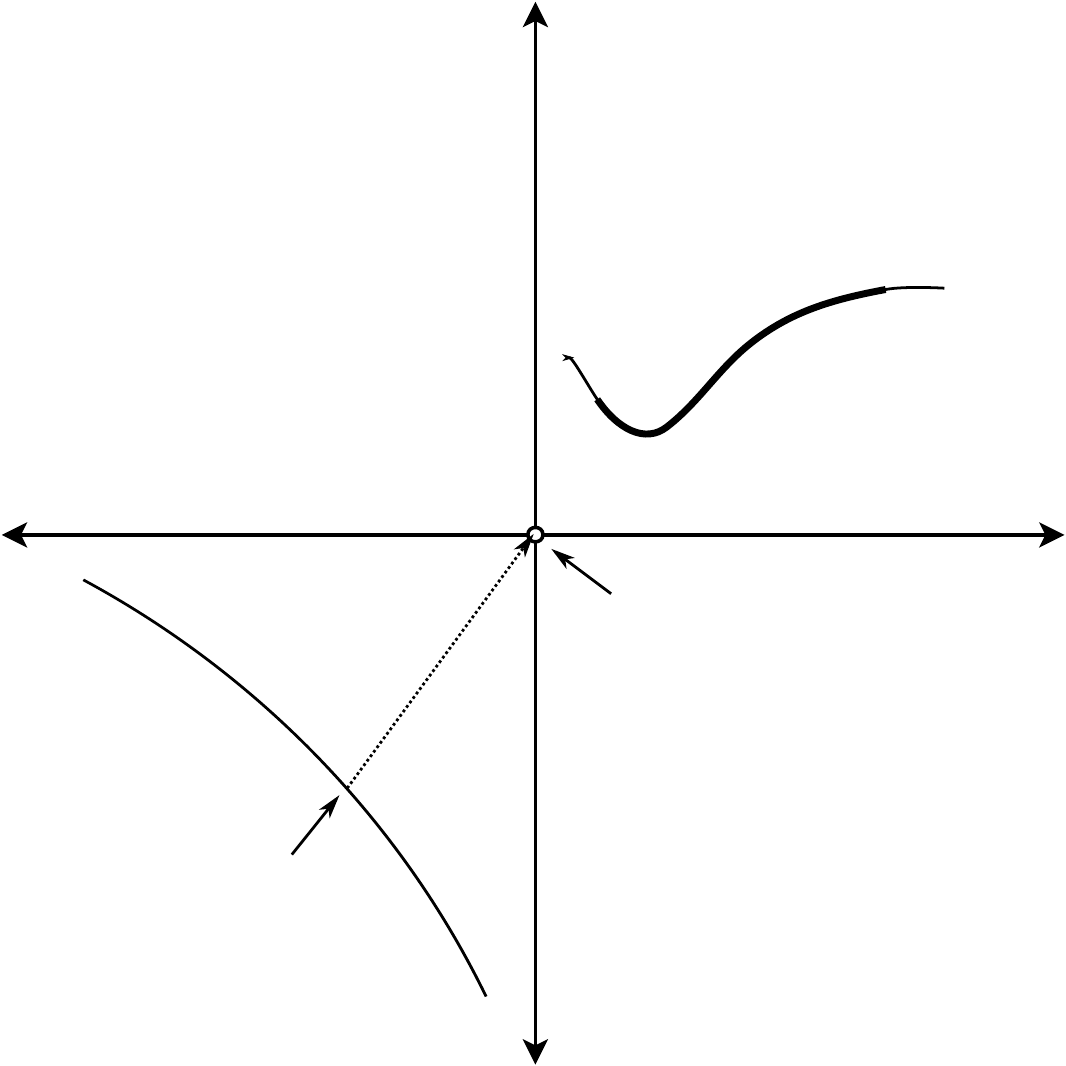}}%
    \put(0.5732899,0.63843648){\color[rgb]{0,0,0}\makebox(0,0)[lb]{\smash{$q_0$}}}%
    \put(0.26058632,0.35179153){\color[rgb]{0,0,0}\makebox(0,0)[lb]{\smash{$q_1$}}}%
    \put(0.58631922,0.40390879){\color[rgb]{0,0,0}\makebox(0,0)[lb]{\smash{$\sqrt{\ph'(\theta)}\; q_1(\tilde\theta)$}}}%
    \put(0.16677525,0.12508141){\color[rgb]{0,0,0}\makebox(0,0)[lb]{\smash{$q_1(\tilde\theta)$}}}%
  \end{picture}%
\endgroup%

\caption{Left side: solution to the finite dimensional minimization problem. In the halfplane below the dotted line, the solution is given by $\sqrt{\ph'(\th)}=0$. On the halfplane
above the dotted line the solution is given by the unique value $\sqrt{\ph'(\th)}$ such that $\sqrt{\ph'(\th)}q_1(\tilde \th)$ lies on the dotted circle. Right side: effect of the reparametrization action on the space of SRVTs. 
}
\label{fig:finite_min}
\end{figure}

Let us look at $\on{dist}_I(c_0, c_1\o\ph)$, first for open curves in the SRV representation. Let $q_i = R(c_i)$. We have $R(c_1 \o \ph) = \sqrt{\ph'}\, q_1 \o \ph$ and
\[
\on{dist}_I(c_0, c_1 \o \ph)^2 = \int_{0}^{2\pi} | q_0(\th) - \sqrt{\ph'(\th)} q_1(\ph(\th)) |^2 \ud \th\;.
\]
Assume that $\ps$ minimizes this expression, fix $\th \in S^1$ and set $\tilde \th = \ps(\th)$. Even though finding $\ps$ has to be done over the whole interval $[0,2\pi]$ simultaneously, it is very instructive to look at just one point at a time. Consider the infimum
\[
\inf_{\sqrt{\ph'(\th)} \geq 0} | q_0(\th) - \sqrt{\ph'(\th)} q_1(\tilde \th) |^2\;.
\]
This is a $d$-dimensional minimization problem, that can be solved explicitly; its solution is visualized in Fig.~\ref{fig:finite_min}. Denote by $\al$ the angle between $q_0(\th)$ and $q_1(\tilde \th)$. If $\frac \pi 2 \leq \al \leq \frac{3\pi}2$, then the infimum is attained for $\sqrt{\ph'(\th)}=0$. In other words, for $q_1(\tilde \th)$ lying in the half-plane ``opposite'' $q_0(\th)$, the optimal reparametrization would scale it to 0.
Next we look at close by points. If the optimal scaling at $\th$ is $\sqrt{\ph'(\th)}=0$ and $\th + \De \th$ is close enough, then the angle between $q_0(\th + \De \th)$ and $q_1(\tilde\th)$ will also lie inside $[\frac \pi 2, \frac{3\pi}2]$ and so $\sqrt{\ph'(\th+\De \th)}=0$ as well. But this would lead to $\ph$ being constant on a whole subinterval of $[0,2\pi]$.

\begin{figure}
\begin{center}
\includegraphics[width=.32\textwidth]{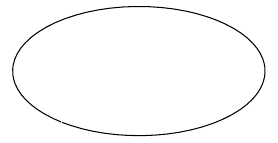}
\includegraphics[width=.32\textwidth]{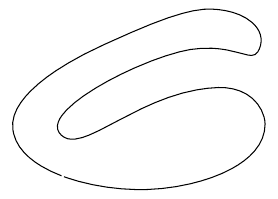}
\includegraphics[width=.32\textwidth]{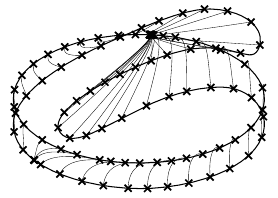}
\end{center}
\caption{Left side: initial curve. Middle figure: target curve. Right figure: minimal geodesic on shape $\mathcal S(S^1,\mathbb R^2)$ between  an ellipse and an ellipse with a large fold. 
One can see, that the fold grows out of a singular point. Image taken from \cite{Bauer2014b}.}
\label{fig:fold}
\end{figure}

The true situation is more complicated than that, in particular for closed curves, where we additionally have the nonlocal constraint $\int_{S^1} q |q| \ud \th =0$ to satisfy. But we do observe the scaling-to-zero behavior in numerical computations; see for example Fig.~\ref{fig:fold}.

\runinhead{Incompleteness}

The key conclusion is this: we should expect the solution $\ps$ of the minimization problem 
$\inf_\ph \on{dist}_I(c_1, c_2 \o \ph)$ to have intervals, where it is constant---that is true, even if we solve the problem on a finite-dimensional approximation space. If $I$ is such an interval and $\ps|_I = \th_I \in S^1$, then this means that the whole segment $c_1(I)$ of the first curve corresponds to the point $c_2(\th_I)$ on the second curve.

Now we can switch $c_1$ and $c_2$. Then the optimal reparametrization is $\ps\inv$. However, since $\ps$ is constant on the interval $I$, its inverse $\ps\inv$ will have a jump at the point $\th_I$. 
What does this mean for minimizing geodesics? If $c(t)$ is a length-minimizing path between $c_2$ and $c_1 \o \ps\inv$, then the point $c_2(\th_I)$ will ``open up'' to the whole segment $c_1(I)$. 

A geodesic is supposed to encode the differences between the shapes represented by $c_1$ and $c_2$ in its initial velocity $\p_t \ga(0)$. However the geodesic starting at $c_2$ sees only the parametrized curve $c_1 \o \ps\inv$ and since $\ps\inv$ has a jump at $\th_I$, jumping over the interval $I$, this interval is missing from the curve $c_1 \o \ps\inv$. How then can the geodesic encode the shape differences, if it does not ``see'' them?

We understand that this is not a rigorous proof\footnote{See \cite[Sect. 4.2]{Michor2008a} for a 
rigorous proof, that this behavior indeed occurs for the metric $G_c(h,k) = \frac{1}{\ell_c} 
\int_{0}^{2\pi} \langle D_s h, D_s k \rangle \ud s$ on the space 
of unparametrized open curves modulo translations.}.    
However there is numerical evidence pointing in the same direction. In Fig.~\ref{fig:fold} we see 
an attempt to numerically compute the minimizing geodesic between an ellipse and an ellipse with a 
large fold, both considered as unparametrized curves. The picture on the right shows the 
point-to-point correspondences after we have minimized over the reparametrization group. We can see 
that indeed one point on the ellipse wants to correspond to the part of the fold, where the tangent 
vectors points in the opposite direction to the tangent vector of the ellipse.     

In the example of Fig.~\ref{fig:fold} we could have cleverly selected the shape with the fold as 
the initial curve and compute our geodesic starting from there. Then there would be no problem of 
one point wanting to become a whole segment. However, for two general curves $c_1$, $c_2$, we would 
expect to encounter mixed behavior: some segments of $c_1$ would collapse to points on $c_2$ and 
other points on $c_1$ would expand to segments of $c_2$.

This effect is not caused by the global ``curvedness'' of the manifold of shapes, it is rather a 
manifestation of the incompleteness.  
We expect to see this as soon as we match two curves whose tangent vectors point in opposite directions. 
Since in the SRV representation distances are measured in the $L^2$-norm, this 
behavior can occur for a pair of curves with arbitrary small (geodesic) distance. 
The fold in Fig.~\ref{fig:fold} can be arbitrary small, but the behavior will be the same. 

Let us fix a representative curve $c$ for the shape $[c]$ and let us look at the space of unparametrized curves through the lens of the exponential map, while standing at the curve $c$. We look at a shape $[c_1]$, represented by a curve $c_1$, by finding a reparametrization $\ps_1$, s.t. $\on{dist}_{I}(c, c_1 \o \ps)$ is minimal, and then we compute $v_1 = \exp_c\inv(c_1 \o \ps_1)$. This vector $v_1$ is what we see, when we look at the shape $[c_1]$. Now, if the reparametrization $\ps$ has jumps, then the curve $c_1 \o \ps_1$ wil miss parts of the shape $[c_1]$; furthermore, there will be several shapes $[c_2]$, distinct from $[c_1]$ only in the part that is missing from $c_1 \o \ps_1$, such that the corresponding optimal representing curves $c_2 \o \ps_2$ coincide with $c_1 \o \ps_1$. This implies that $v_2 = \exp_c\inv(c_2 \o \ps_2)$ coincides with $v_1$, while the shapes $[c_2]$ and $[c_1]$ differ. In other words we look at different shapes, but see the same thing. In fact there are many regions in shape space, that cannot be distinguished using the exponential map and these regions start arbitrary close to the starting shape $[c]$.

\runinhead{Joint Reparametrizations}
It is possible that searching for \emph{one} reparametrization is the wrong problem. Mathematically an equivalent way to define the geodesic distance on $B(S^1,\R^d)/\on{Tra}$ is via
\[
\on{dist}_B(C_1, C_2) = \inf_{\ph_1,\ph_2 \in \on{Diff}(S^1)} \on{dist}_I(c_1 \o \ph_1, c_2 \o \ph_2)\;.
\]
Using the invariance of $\on{dist}_I$ under reparametrizations we can recover~\eqref{eq:shape_distance}. Now we are looking for reparametrizations of both curves, such that the infimum is attained. The advantage of this approach is that we can avoid jumps, that would necessarily appear in the one--reparametrization strategy, by instead setting the other reparametrization to be constant on the corresponding interval. 

The underlying behavior does not change: points on one curve can be matched to intervals on the other and vice versa. If $\ps_1,\ps_2$ represent a pair of optimal reparametrizations for two curves $c_1,c_2$ and $\ps_1$ is constant on the interval $I$ with $\ps_1|_I = \th_I$, then the point $c_1(\th_I)$ will correspond to the interval $c_2(\ps_2(I))$. Instead of jumping over the interval $\ps_2(I)$, we now reparametrize $c_1$ and the reparametrized curve $c_1\o\ps_1$ waits at $\th_I$ until $\ps_2$ has moved past $I$.

The strategy of joint reparametrizations was proposed in \cite{Robinson2012, Lahiri2015_preprint}, where the authors consider only curves without the periodicity constraint. Even for open curves, it is not known, whether for any two absolutely continuous curves, there exists a pair of reparametrizations realizing the geodesic distance; in \cite{Lahiri2015_preprint} this is shown only under the additional assumption, that one of the curves is piecewise linear.

\runinhead{Other $H^1$-metrics}
There are many different $H^1$-metrics to choose from. For a start, the constants $1$ and $\frac 14$ are rather arbitrary and we could look at the full family of metrics of the form
\[
G_c(h, k) = \int_{S^1} a^2 \langle D_s h^\perp, D_s k^\perp \rangle 
+ b^2 \langle D_s h,v\rangle \langle D_s k,v\rangle \ud s\;,
\]
with $a, b > 0$. This family has been given the name \emph{elastic metrics} and has been studied for plane curves in \cite{Mio2007, Bauer2014b}. All metrics in this family are uniformly equivalent, i.e., if $G$ and $H$ are two elastic metrics with possibly different constants $a,b$, there exists a constant $C$, such that
\[
C\inv H_c(h,h) \leq G_c(h,h) \leq C H_c(h,h)\;,
\]
holds for all curves $c$ and all tangent vectors $h$.

All $H^1$-metrics can be made invariant with respect to scalings by multiplying them with an appropriate power of the length $\ell(c)$; for example the following metric is scale-invariant,
\[
G_c(h, k) = \int_{S^1} \ell_c^{-3} \langle h, k \rangle 
+ \ell_c\inv \langle D_s h, D_s k \rangle \ud s\;.
\]
A modified $H^1$-metric was introduced in \cite{Sundaramoorthi2011} with the property that scalings, translations and general deformations of the curve are all orthogonal to each other; this metric was then applied to tracking moving objects.

We have looked at completeness properties only for the $H^1$-metric corresponding to the SRVT; a similar, more rigorous, discussion can be found for the scale-invariant version of elastic metric corresponding to the choice $a=b$ in \cite{Michor2008a} and we conjecture, that all $H^1$-metrics share the same qualitative behavior. We will see in the next section, what happens, if we additionally penalize second and higher derivatives of the tangent vectors.

\section{Higher Order Sobolev Metrics}\label{sec:higherorder}

Riemannian metrics involving first derivatives of the tangent vectors can lead to very efficient computations, but some of their mathematical properties are less convenient. Now we will make the metric dependent on higher derivatives. It is not easy to give a definition of a general Sobolev-type Riemannian metric, that is both general enough to encompass all known examples and concrete enough so that we can work easily with it. We will approach this class by looking at families of examples instead, noting common features as well as differences.

A very useful family is that of {\itshape Sobolev metrics with constant coefficients}. These are metrics of the form
\begin{equation}\label{eq:constant_coeff}
G_c(h,k) = \int_{S^1} a_0 \langle h, k \rangle + a_1 \langle D_s h, D_s k\rangle + \dots +
a_n \langle D_s^n h, D_s^n k \rangle \ud s\;,
\end{equation}
with constants $a_0,\dots,a_n$. 
The largest $n$, such that $a_n \neq 0$ is called the {\itshape order} of the metric. We require $a_j \geq 0$ for the metric to be positive semi-definite, $a_n > 0$ for it to be a metric of order $n$ and $a_0 > 0$ for it to be non-degenerate. If $a_0=0$, then constant tangent vectors are in the kernel of $G$ and, provided there is at least one non-zero coefficient, $G$ defines a non-degenerate metric on the quotient space $\on{Imm}(S^1,\R^d)/\on{Tra}$ of regular curves modulo translations. Most of the metrics encountered in Sect.~\ref{sec:SRVT} were of this type.

Using integration by parts we can rewrite \eqref{eq:constant_coeff} to obtain
\begin{equation}\label{eq:constant_coeff2}
G_c(h,k) = \int_{S^1} a_0 \langle h,  k \rangle + a_1 \langle -D^2_s h, k\rangle + \dots +
a_n \langle (-1)^n D_s^{2n} h,  k \rangle \ud s\;,
\end{equation}
enabling us to write the metric in the form $G_c(h,k)=\int_{S^1} \langle L_ch,k\rangle \ud s$, with $L_c=\sum_{j=0}^n (-1)^j a_j\, D_s^{2j}$, a differential operator of order $2n$.

\runinhead{Metrics with Nonconstant Coefficients} 
We could loosen our restrictions on the coefficients $a_j$ and permit them to be functions, that depend on the curve $c$ and quantities derived from it, e.g., $a_j = a_j(\ell_c, \ka, D_s \ka, \dots)$. In Sect.~\ref{sec:L2metric} we have 
considered such metrics of order zero, the almost local metrics.
Several examples of higher order metrics with nonconstant coefficients have been investigated in the literature. 
The completeness properties of first and second order metrics with coefficients depending on the length are studied in \cite{Mennucci2008}. The idea in \cite{Shah2013} and \cite{Bauer2014_preprint} is to decompose a tangent vector $h=h^\parallel + h^\perp$ into a part tangent and a part normal to the curve and consider derivatives of these quantities. Some special examples of second order metrics can be found in \cite{Bauer2014c}.

\runinhead{Geodesic Equation} The geodesic equation of a Sobolev metric with constant coefficients is a nonlinear PDE, second order in time and of order $2n$ in the space variable. It is given by
\begin{equation*}
\begin{split}
\p_t \left(\sum_{j=0}^n (-1)^j a_j \,|c'|\, D_s^{2j} c_t\right) &=
-\frac{a_0}2 \,|c'|\, D_s\left( \langle c_t, c_t \rangle D_s c \right) \\
&\qquad{} + \sum_{k=1}^n \sum_{j=1}^{2k-1} (-1)^{k+j} \frac{a_k}{2}\, |c'|\, D_s
\left(\langle D_s^{2k-j} c_t, D_s^j c_t \rangle D_s c \right)\;.
\end{split}
\end{equation*}
We can see that if $a_j=0$ for $j \geq 1$, then this equation reduces to the geodesic equation \eqref{eq:geode_L2} of the $L^2$-metric. The left hand side of the geodesic equation is the time derivative of the momentum, $L_c c_t\,|c'|$. For metrics of order $n\geq 1$ the geodesic equation is locally well-posed \cite{Michor2007}.

Now we come to the main difference between Sobolev metrics of order one and metrics of higher order. In a nutshell, first order Sobolev metrics are only \emph{weak} Riemannian metrics, while Sobolev metrics of higher order, if extended to a suitable, larger space, are \emph{strong} Riemannian metrics.

\runinhead{Weak Sobolev Metrics}
Let $G$ be a Sobolev metric of order one,
\[
G_c(h,k) = \int_{S^1} \langle h, k \rangle + \langle D_s h, D_s k \rangle \ud s
= \int_{S^1} \langle h, k \rangle |c'|+ \langle h', k' \rangle |c'|\inv \ud \th\;.
\]
Fix the curve $c$ and look at the inner product $G_c(\cdot,\cdot)$. The natural space to define $G_c(\cdot,\cdot)$ is the Sobolev space $H^1(S^1,\R^d)$ of functions with square-integrable derivatives, together  with the inner product
\[
\langle h, k \rangle_{H^1} = \int_{S^1} \langle h, k \rangle + \langle h', k' \rangle \ud \th\;.
\]
If $c$ is smooth enough, say $c \in C^1$, then we see that $G_c(\cdot,\cdot)$ defines an inner product on $H^1$, which is equivalent to the standard inner product. Unfortunately we cannot allow $c$ itself to be an $H^1$-function. We need uniform control on the derivative $c'$ to guarantee that the integral $\int_{S^1} \langle h',k'\rangle |c'|\inv \ud \th$ is finite, but the $H^1$-norm does not provide that. The best we can do is to extend $G$ to the space 
\[
G: C^1\!\on{Imm}(S^1,\R^d) \x H^1(S^1,\R^d) \x H^1(S^1,\R^d) \to \R\;.
\]
In this sense $G$ is a {\itshape weak Riemannian metric}\footnote{An infinite-dimensional Riemannian manifold $(M,g)$ is called \emph{strong}, if $g$ induces the natural topology on each tangent space or equivalently, if the map $g: TM \to (TM)'$ is an isomorphism. If $g$ is merely a smoothly varying nondegenerate bilinear form on $TM$ we call $(M,g)$ a \emph{weak} Riemannian manifold, indicating that the topology induced by $g$ can be weaker than the natural topology on $TM$ or equivalently $g:TM \to (TM)'$ is only injective.}; the topology induced by the inner product $G_c(\cdot,\cdot)$, in this case the $H^1$-topology, is weaker than the manifold topology, here the $C^\infty$- or $C^1$-topology.

\runinhead{Strong Sobolev Metrics} The situation is different for Sobolev metrics with constant coefficients of order 2 or higher. Let us look at an example:
\[
G_c(h,k) = \int_{S^1} \langle h, k \rangle + \langle D_s^2 h, D_s^2 k \rangle \ud s\;,
\]
with $D_sh = \frac{1}{|c'|} h'$ and $D_s^2 h = \frac{1}{|c'|} h'' - \frac{\langle c',c''\rangle}{|c'|} h'$. Again the natural space for $G_c(\cdot,\cdot)$ is the Sobolev space $H^2(S^1,\R^d)$ and it would appear that we need $c \in C^2$ for the inner product to be well-defined. However a careful application of Sobolev embedding and multiplier theorems---see~\cite[Sect.~3.2]{Bruveris2014}---shows that we can extend $G$ to a smooth inner product on the space
\[
G : \mc I^2(S^1,\R^d) \x H^2(S^1,\R^d) \x H^2(S^1,\R^d) \to \R\;;
\]
here we denote by $\mathcal I^2(S^1,\R^d) = \{ c \in H^2\,:\, c'(\th) \neq 0\; \forall \th \in S^1\}$ the space of $H^2$-curves with non-vanishing tangent vectors. The crucial fact is the Sobolev embedding $H^2 \hookrightarrow C^1$, implying that the $H^2$-norm controls first derivatives uniformly. This also implies that $\mc I^2$ is an open set in $H^2$. Thus $G$ becomes a {\itshape strong Riemannian metric} on $\mc I^2(S^1,\R^d)$; the topology induced by each inner product $G_c(\cdot,\cdot)$ coincides with the manifold topology.

Similarly Sobolev metrics of order $n$ with constant coefficients induce strong metrics on the space $\mc I^n(S^1,\R^d)$ of regular Sobolev curves, provided $n \geq 2$.

Note however that Sobolev metrics of order 2 and higher are strong metrics only when considered on the larger space $\mc I^n(S^1,\R^d)$, not on the space $\on{Imm}(S^1,\R^d)$ of smooth curves. On $\on{Imm}(S^1,\R^d)$ the metric is still a weak metric. That said the difference between metrics of order 1 and those of higher order is that for higher order metrics we are able to pass to the larger space $\mc I^n(S^1,\R^d)$---one could say we are ``completing'' the space of smooth curves---on which it becomes a strong metric, while for first order metrics such a completion does not exist.

\runinhead{Properties of Strong Metrics} The following is a list of properties we get ``for free'' simply by working with a smooth, strong Riemannian metric as opposed to a weak one:
\begin{itemize}
\item
The Levi-Civita covariant derivative exists and the geodesic equation has local solutions, which depend smoothly on the initial conditions.
\item
The exponential map exists and is a local diffeomorphism.
\item
The induced geodesic distance is point-separating
and generates the same topology as the underlying manifold.
\end{itemize}

The theory of strong, infinite-dimensional Riemannian manifolds is described in \cite{Lang1999} and \cite{Klingenberg1995}. For weak Riemannian manifolds all these properties have to be established by hand and there are examples, where they fail to hold. The geodesic distance for the $L^2$-metric discussed in Sect.~\ref{sec:L2metric} on the space of curves vanishes identically \cite{Michor2006c, Michor2005}, it is not known whether the geodesic equation for the $L^2$-metric has solutions and \cite{Bauer2014d} presents a weak Riemannian manifold that has no Levi-Civita covariant derivative.

We would like to note that the distinction between weak and strong Riemannian metrics arises only 
in infinite-dimensions. Every Riemannian metric on a finite-dimensional manifold is strong. 
Therefore phenomena, like the vanishing of the geodesic distance or the failure of geodesics to 
exist, can only arise in infinite dimensions. For better or for worse, this is the setting, where 
the joy and pain of shape analysis occurs.

\runinhead{Completeness Properties} What are the operations on the manifold of curves, that are used in shape analysis?
\begin{enumerate}[(1)]
\item
We want to use the Riemannian exponential map $\on{exp}_c : T_c M \to M$ to pass between the tangent space at one point and the manifold itself. Its inverse $\on{log}_c = \on{exp}_c \inv$ allows us to represent the nonlinear manifold or at least a part thereof in a vector space.
\item
Given two curves, we want to compute the geodesic distance between them, that is, the length of the shortest path joining them. Often we are also interested in the shortest path itself; it can be used to transfer information between two curves and the midpoint can serve as the average of its endpoints.
\item
Given a finite set $\{c_1,\dots,c_n\}$ of curves, we are interested in the average of this set. On a Riemannian manifold this is usually the Fr\'echet or Karcher mean, i.e., a curve $c^\ast$, minimizing the sum of squared distances, $c^\ast = \on{argmin}_c \sum_i d(c,c_i)^2$, where $d$ is the geodesic distance.
\end{enumerate}

This is by no means an exhaustive list, but describes rather the basic operations, one wants to perform. The ability to do so places conditions on the manifold of curves and the Riemannian metric. Let us look at these conditions.
\begin{enumerate}[(1)]
\item
\label{geod_compl}
On a general Riemannian manifold the exponential map $\on{exp}_c : U \to M$ is defined only on an open neighborhood $U \subseteq T_c M$ of 0 and it is rarely known exactly how $U$ looks like. For us to be able to freely map tangent vectors to curves, we want a globally defined exponential map. Since the exponential map is defined as $\on{exp}_c(h) = \ga(1)$, where $\ga$ is the geodesic with initial conditions $(c,h)$ and using the property $\on{exp}_c(th) = \ga(t)$, we can see that a globally defined exponential map is equivalent to requiring that geodesics exist for all time. This property is called {\itshape geodesic completeness}.
\item[(1')]
Asking for the exponential map to be invertible is more difficult. On a strong Riemannian manifold this is always the case locally. Imposing it globally is a very restrictive condition. The Weil--Peterson metric \cite{Mumford2006} comes closest; it is a metric with negative sectional curvature, meaning that the derivative of the exponential map is everywhere invertible.
\item
\label{ex_min_geod}
Here we want to know, whether any two curves can be joined by a minimizing geodesic, i.e., a geodesic, whose length realizes the geodesic distance. On a finite-dimensional manifold geodesic completeness would imply this property; this is not the case in infinite dimensions \cite{Atkin1975}. This is not to say, that we cannot hope for minimizing geodesics to exist, but rather, that it will have to be proven independently of (\ref{geod_compl}).
\item
\label{metr_compl}
Ensuring that the Fr\'echet mean exists for all finite collections of curves is difficult. But there is a theorem \cite{Azagra2005} stating that the mean exists and is unique on a dense subset of the $n$-fold product $M \x \dots \x M$, provided the manifold is {\itshape metrically complete}. This means that the manifold $(M,d)$ together with the induced geodesic distance is complete as a metric space.
\end{enumerate}

The properties (\ref{geod_compl}), (\ref{ex_min_geod}) and (\ref{metr_compl}) for Riemannian manifolds are called completeness properties. In finite dimensions the theorem of Hopf--Rinow states that (\ref{geod_compl}) and (\ref{metr_compl}) are equivalent and either one implies (\ref{ex_min_geod}). For infinite-dimensional strong Riemannian manifolds the only remaining implication is that metric completeness implies geodesic completeness.\footnote{A counterexample, showing that in infinite dimensions metric and geodesic completeness together do not imply existence of minimizing geodesics can be found in \cite{Atkin1975}; similarly, that geodesic completeness does not imply metric completeness is shown in \cite{Atkin1997}}.

\runinhead{Completeness for Sobolev Metrics}
Let us look at the situation for Sobolev metrics on the space of parametrized curves. We have argued in Sect.~\ref{sec:SRVT} that we shouldn't expect $H^1$-metrics to be geodesically or metrically complete. Things look better for Sobolev metrics of higher order. In fact it is shown in \cite{Bruveris2014} and \cite{Bruveris2014b_preprint} that these metrics satisfy all the above mentioned completeness properties. 
\footnote{A related result holds for the space of curves of bounded second variation, together with a Finsler $BV^2$-metric. It is shown in \cite{Vialard2014_preprint} that any two curves in the same connected component of the space $BV^2(S^1,\R^2)$ can be joined by a length-minimizing path.}

\begin{theorem}
Let $n \geq 2$ and let $G$ be a Sobolev metric of order $n$ with constant coefficients. Then
\begin{enumerate}[(1)]
\item
$(\mc I^n(S^1,\R^d), G)$ is geodesically complete;
\item
Any two elements in the same connected component of $\mc I^n(S^1,\R^d)$ can be joined by a minimizing geodesic;
\item
$(\mc I^n(S^1,\R^d), \on{dist}_{\mc I})$ is a complete metric space.
\end{enumerate}
\end{theorem}

The geodesic equation for Sobolev metrics has a smoothness preserving property. If the initial conditions are smoother that $H^n$---let us say the initial curve and initial velocity are $C^\infty$---then the whole geodesic will be as smooth as the initial conditions. Therefore the space $(\on{Imm}(S^1,\R^d),G)$ of smooth immersions with a Sobolev metric of order $n$ is also geodesically complete.

\runinhead{Completeness for Unparametrized Curves} Similar completeness properties hold for unparametrized curves. The correct space, where to look for completeness, is the quotient
\[
\mc B^n(S^1,\R^d) = \mc I^n(S^1,\R^d) / \mc D^n(S^1)\;,
\]
of regular curves of Sobolev class $H^n$ modulo reparametrizations of the same regularity; the group $\mc D^n(S^1) = \{ \ph \in H^n(S^1,S^1) \,:\, \ph'(\th) > 0\}$ is the group of $H^n$-diffeomorphisms. We have the following theorem \cite{Bruveris2014b_preprint}:

\begin{theorem}
\label{shape_compl}
Let $n \geq 2$ and let $G$ be a Sobolev metric of order $n$ with constant coefficients. Then
\begin{enumerate}[(1)]
\item
$(\mc B^n(S^1,\R^d), \on{dist}_{\mc B})$ with the quotient metric induced by the geodesic distance on $(\mc I^n, G)$ is a complete metric space;
\item
\label{ex_reparam}
Given $C_1, C_2 \in \mc B^n$ in the same connected component, there exist $c_1, c_2 \in \mc I^n$ with $c_1 \in \pi\inv(C_1)$ and $c_2 \in \pi\inv(C_2)$, such that
\[
\on{dist}_{\mc B}(C_1, C_2) = \on{dist}_{\mc I}(c_1, c_2)\;;
\]
equivalently, the infimum in
\[
\on{dist}_{\mc B}(\pi(c_1), \pi(c_2)) = \inf_{\ph \in \mc D^n(S^1)} \on{dist}_{\mc I}(c_1, c_2 \o \ph)
\]
is attained;
\item
Any two shapes in the same connected component of $\mc B^n(S^1,\R^d)$ can be joined by a minimizing geodesic.
\end{enumerate}
\end{theorem}

The only drawback is that the space $\mc B^n(S^1,\R^d)$ of Sobolev shapes is not a manifold any 
more\footnote{This has to do with smoothness properties of the composition map in Sobolev spaces. 
While the smooth reparametrization group $\on{Diff}(S^1)$ acts smoothly on the space 
$\on{Imm}(S^1,\R^d)$ of smooth curves \cite{KrieglMichor97}, the group of Sobolev reparametrizations $\mc D^n(S^1)$ acts 
only continuously on Sobolev curves $\mc I^n(S^1,\R^d)$. Moreover, the smooth shape space $B(S^1,\R^d)$ 
is (apart from isolated singularities) a manifold, but the Sobolev shape space $\mc B^n(S^1,\R^d)$ 
is only a topological space. This is the price we have to pay for completeness. See \cite[Sect.~6]{Bruveris2014b_preprint} for details.}. It is however a topological space and equipped with the geodesic distance function a metric space. We have to understand a minimizing geodesic in the sense of metric spaces, i.e., a curve $\ga : I \to \mc B^n$ is a minimizing geodesic, if     
\[
\on{dist}_{\mc B}(\ga(t), \ga(s)) = \la |t-s|
\]
holds for some $\la > 0$ and all $t,s \in I$.

We would like to point out that part (\ref{ex_reparam}) of Thm.~\ref{shape_compl} is the counterpart of the incompleteness discussion in Sect.~\ref{sec:SRVT}. This theorem states that given two shapes, represented by two parametrized curves, we can find an optimal reparametrization $\ph$ of the second curve. The fact that $\ph \in \mc D^n(S^1)$ guarantees that $\ph$ is at least a $C^1$-diffeomorphism of the circle; thus no intervals are collapsed to single points and neither $\ph$ nor $\ph\inv$ has any jumps.

\runinhead{Numerical Computations}
We have argued that Sobolev metrics of sufficiently high order have nice mathematical properties, which are relevant to applications in shape analysis. What we have not done is present convincing applications, showcasing their superior performance. This is because the numerical computation of minimizing geodesics between parametrized or unparametrized curves is still an open problem. First attempts at discretizing special cases of $H^2$-metrics can be found in \cite{Bauer2014_preprint,Vialard2014_preprint}. While first order metrics have nice representations in terms of the SRVT or the basic mapping \cite{Michor2008a}, which greatly simplifies the numerics, there is no such analogon for higher order metrics\footnote{In \cite{Bauer2014c} a representation of second order metrics, similar to the SRVT was developed. However, image of the resulting transformations have infinite co-dimension, which, compared to the SRVT, complicates the situation.} Finding a robust and stable discretization of metrics of order 2 and higher remains a challenge.

\section{Riemannian Metrics on the Space of Surfaces}

In the previous sections we have presented several reparametrization invariant metrics on the space of curves. We want to conclude the exposition with a short excursion to the space of regularly parametrized surfaces, i.e., $\on{Imm}(M,\R^d) = \{ f \in C^\infty(M,\R^d) \,:\, T_xf\text{ injective } \forall x \in M\}$ with $M$ a compact 2-dimensional manifold without boundary. Typical choices for $M$ are the sphere $S^2$ and the torus $S^1 \x S^1$.
In particular we will describe how the previously described metrics can be generalized from $\on{Imm}(S^1,\R^2)$ to $\on{Imm}(M,\R^d)$.
We will follow the presentation of \cite{Bauer2011b,Bauer2012a}.

\runinhead{Sobolev metrics}
To generalize Sobolev metrics~\eqref{eq:constant_coeff2} from the space of curves to the space of surfaces, we need the right replacements for the arc length derivative $D_s$ and integration $\ud s$. For an immersion $f \in \on{Imm}(M,\R^d)$, we denote by $g = g^f = f^\ast\langle \cdot, \cdot \rangle$ the pullback of the Euclidean metric to $M$. This makes $(M,g)$ into a Riemannian manifold with a Laplace operator $\De^g = -\on{div}^g \o \on{grad}^g$ and volume form $\on{vol}^g$. We will use $\De^g$ and $\on{vol}^g$ as replacements for $-D_s^2$ and $\ud s$. A Sobolev metric of order $n$ on $\on{Imm}(M,\R^d)$ is given by
\begin{equation}\label{eq:sobmetric_surfaces}
G_f(h,k)= \int_{M} a_0 \langle h, k \rangle + a_1 \langle \Delta^g h, k\rangle + \dots +
a_n \langle (\Delta^g)^n h, k \rangle \on{vol}^g\;;
\end{equation}
here the tangent vectors $h,k$ are seen as maps $h,k : M \to \R^d$ and the Laplace operator acts on each coordinate separately\footnote{
On $\on{Imm}(M,\R^d)$ there is no canonical way to define Sobolev metrics. We could have also used the definition
\[
G_f(h,k) = \int_M a_0 \langle h,k \rangle + \dots + a_n \langle (\nabla^g)^n h, (\nabla^g)^n k \rangle \on{vol}^g\;,
\]
where $\nabla^g$ is the covariant derivative on $M$. Then the differential operator associated to this metric is
$L_f = \sum_{j=0}^n a_j  \left((\nabla^g)^*\right)^n \left(\nabla^g\right)^n\otimes \on{vol}^g$, with $(\nabla^g)^*=\on{Tr}(g^{-1}\nabla^g)$ denoting the adjoint of $\nabla^g$. When $n \geq 2$, the operator $\left((\nabla^g)^*\right)^n \left(\nabla^g\right)^n$ differs
from $(\Delta^g)^n$ in the lower order terms---they are connected via the Weitzenb\"ock formulas.}.

The associated  operator $L_f$, allowing us to write $G_f(h,k)= \int_{S^1}\langle L_f
h,k \rangle\ud s$, is given by $L_f = \sum_{j=0}^n a_j \left(\Delta^g\right)^j \otimes \on{vol}^g$. Every metric in this family 
is invariant under the action of the diffeomorphism group $\on{Diff}(M)$ and induces a Riemannian metric on the quotient space of unparametrized surfaces
$B(M,\R^3)=\on{Imm}(M,\R^3)/\on{Diff}(M)$.

Similarly as in the previous section one can also allow the coefficients $a_j$ to depend on the surface $f$. Coefficients depending on the total volume, the mean curvature and the Gau{\ss} curvature have been considered 
in \cite{Bauer2012d}. The class of almost local metrics on surfaces has also been studied \cite{Bauer2012a}. 

\runinhead{Geodesic Distance}
For the geodesic distance we obtain similar results as for curves:
the geodesic distance vanishes for the $L^2$-metric on $\on{Imm}(M,\R^d)$. Both higher order Sobolev metrics and almost local metrics depending on mean curvature or total volume induce point-separating distances on the space of unparametrized surfaces. Whether higher order Sobolev metrics induce a point-separating distance function on $\on{Imm}(M,\R^d)$ itself is not known.

\runinhead{Geodesic Equation}
The formulas for the geodesic equations for Sobolev metrics \eqref{eq:sobmetric_surfaces} become very quickly very technical; see \cite{Bauer2011b}. As an example we 
present the geodesic equation of the $H^1$-metric with coefficients $a_0=a_1=1$; this is the metric induced by the operator field $L_f = 1+\Delta^g$.
$$\begin{aligned}
\p_t \left(L_f f_t \otimes \on{vol}^g\right) &=\bigg(
\on{Tr}\big(g^{-1} S^f g^{-1} \langle \nabla^g f_t,\nabla^g f_t \rangle  \big)
-\frac{1}{2}\on{Tr}\big(g^{-1} \nabla^g\langle \nabla^g f_t ,f_t\rangle  \big).H^f\\&\quad
-\frac12\langle L_f f_t,f_t\rangle. H^f
-Tf.\langle L_f f_t,\nabla^g f_t\rangle^\sharp\bigg) \otimes \on{vol}^g\;;
\end{aligned}$$
here $S^f$ denotes the second fundamental form, $H^f=\on{Tr}(g^{-1} S^f)$ the vector valued mean curvature and $\nabla^g$ the covariant derivative of the surface $f$.

\runinhead{Outlook}
On spaces of surfaces the theory of Sobolev metrics is significantly less developed than on the space of curves. 
We conjecture that Sobolev metric of order $n \geq 3$ will again be strong Riemannian metrics on the space $\mc I^n(M,\R^d)$ of Sobolev surfaces. 
Nothing is known about completeness properties of these metrics.

An analogue of the SRVT transform has been developed for surfaces in \cite{Jermyn2012,Kurtek2010}. 
However questions regarding the invertibility of the transform and the characterization of its image remain open. 
So far no numerical experiments for higher order Sobolev metrics on the space of surfaces have been conducted.

\begin{acknowledgement}
This work was supported by the Erwin Schr\"odinger Institute. Martin Bauer was supported by the `Fonds zur F\"orderung der wissenschaftlichen Forschung, Projekt P 24625'.
\end{acknowledgement}


\end{document}